\documentclass[11pt]{amsart}
\usepackage{amssymb}
\usepackage{amsmath}

\newcommand{\bbC}{{\mathbb C}}
\newcommand{\bbZ}{{\mathbb Z}}

\newcommand{\bbP}{{\mathbb P}}

\newcommand{\bbQ}{{\mathbb Q}}
\newcommand{\Oh}{{\mathcal O}}

\DeclareMathOperator{\HH}{H}

\DeclareMathOperator{\JJ}{J}

\DeclareMathOperator{\VV}{V}
\DeclareMathOperator{\UU}{U}
\DeclareMathOperator{\CH}{CH}
\DeclareMathOperator{\Pic}{Pic}

\DeclareMathOperator{\Ext}{Ext}
\DeclareMathOperator{\Hom}{Hom}

\DeclareMathOperator{\im}{Image}

\DeclareMathOperator{\reg}{reg}

\newcommand{\onto}{\twoheadrightarrow}

\newcommand{\into}{\hookrightarrow}
\newcommand{\by}[1]{\xrightarrow{#1}}
\newcommand{\ext}[1]{\operatorname{\stackrel{#1}{\wedge}}}
\newcommand{\tensor}{\otimes}
\newcommand{\isom}{\cong}

\newcommand{\sE}{{\mathcal E}}

\newcommand{\sH}{{\mathcal H}}
\newcommand{\sHom}{\mbox{${\sH}{om}$}}
\newcommand{\sExt}{\mbox{${\sE}{xt}$}}

\theoremstyle{plain}
\newtheorem{lemma}{Lemma}
\newtheorem{prop}{Proposition}
\newtheorem{thm}{Theorem}
\newtheorem{conj}{Conjecture}
\newtheorem{cor}{Corollary}

\newtheorem*{thma}{Theorem A}
\newtheorem*{thmb}{Theorem B}
\newtheorem*{corc}{Corollary C}

\theoremstyle{remark}
\newtheorem{remark}{Remark}

\newenvironment{diagram}[1]{\arraycolsep=\doublerulesep\begin{array}{#1}
    }{\end{array}}

\long\def\comment#1{}

\begin{document}
\title[Arithmetically Cohen-Macaulay Bundles]{Arithmetically
Cohen-Macaulay Bundles on complete intersection varieties of
sufficiently high multidegree} 
\author{Jishnu Biswas}
\address{Theoretical Statistics and Mathematics Unit,
Indian Statistical Institute, 8th Mile Mysore Road,
Bangalore 560 059, INDIA.}
\email{jishnu@isibang.ac.in}

\author{G. V. Ravindra}
\address{Department of Mathematics, Indian Institute of Science,
Bangalore 560 012, INDIA.}
\email{ravindra@math.iisc.ernet.in}
\address{Current Address:
Department of Mathematics and Computer Science,
1 University Boulevard, University of Missouri,
St. Louis, MO 63121, USA.}


\begin{abstract}
Recently it has been proved that any arithmetically Cohen-Macaulay
(ACM) bundle of rank two on a general, smooth hypersurface of degree
at least three and dimension at least four is a sum of line
bundles. When the dimension of the hypersurface is three, a similar
result is true provided the degree of the hypersurface is at least
six. We extend these results to complete intersection subvarieties by
proving that any ACM bundle of rank two on a general, smooth complete
intersection subvariety of sufficiently high multi-degree and
dimension at least four splits. We also obtain partial results in the
case of threefolds.
\end{abstract}
\maketitle

\section{Introduction}
We work over the field of complex numbers which shall be denoted by
$\bbC$.

The motivation for the results of this article lie in the study of
certain conjectures of Griffiths and Harris on the structure of curves
in hypersurfaces in $\bbP^4$. These conjectures can be viewed as a
generalisation of the Noether-Lefschetz theorem which we recall now.

\begin{thm}[Noether-Lefschetz theorem]\label{nlt}
Let $X\subset\bbP^3$ be a smooth, very general hypersurface of degree
$d\geq 4$. Then any curve $C\subset X$ is a complete intersection,
i.e., $C=X\cap{S}$ where $S\subset \bbP^3$ is a hypersurface.
\end{thm}

Inspired by the above theorem, Griffiths and Harris (see \cite{GH})
made a series of conjectures in decreasing strength about the
structure of $1$-cycles in $X$, the strongest one of which is the
following.

\begin{conj}\label{ghconj}
Let $X\subset \bbP^4$ be a general, smooth hypersurface of degree
$d\geq 6$, and $C\subset X$ be any curve. Then $C=X\cap{S}$ where
$S\subset \bbP^4$ is a surface.
\end{conj}

This conjecture was proved to be false by Voisin (see \cite{V}). In
fact, she showed that pursuing a certain line of thought, which we
describe below, weaker versions of this conjecture are also false.

Notice that in the Noether-Lefschetz situation, for a smooth curve $C\subset
X$, the normal bundle sequence 
\begin{eqnarray}\label{nbs}
0 \to N_{C/X} \to N_{C/\bbP^3} \to \Oh_C(d) \to 0
\end{eqnarray}
splits. Griffiths and Harris investigated the splitting of this normal bundle sequence and
proved (see \cite{GH1}) the following characterisation 

\begin{thm}\label{nbs1}
Let $X \subset \bbP^3$ be a smooth hypersurface of degree $d$ and let
$C\subset X$ be a smooth curve. Then the normal bundle sequence \eqref{nbs}
splits if and only if $C\subset X$ is a complete intersection.
\end{thm}

Unfortunately, the situation is not as simple in higher
dimensions. Let $X\subset \bbP^4$ be a smooth hypersurface of degree
$d$ and $C\subset X$ be any smooth curve. It is not hard to see that
if $C=X\cap{S}$ where $S\subset\bbP^4$ is a surface, then the normal
bundle sequence for the inclusion $C\subset X \subset \bbP^{4}$ splits.
Furthermore if $X_2$ denotes the first order thickening of $X$ in
$\bbP^4$, then the splitting of the above normal bundle sequence
implies the splitting of the sequence
$$0\to N_{C/X} \to N_{C/X_2} \to N_{X/X_2}|_C \to 0.$$ 
By Lemma 1 in \cite{MPRV}, the splitting of this normal bundle
sequence implies that there exists $D\subset X_2$, a one dimensional
subscheme ``extending'' $C$ i.e., $C=D\cap{X}$.

It is this weaker splitting that Voisin investigates.  In \cite{V},
she proves the following
\begin{prop}
Let $X\subset \bbP^4$ be a smooth hypersurface of degree $d>1$. There
exist smooth curves $C\subset X$ such that the normal bundle sequence
for the inclusions $C\subset X \subset X_2$ (and hence for the
inclusions $C\subset X\subset \bbP^4$) does not split. In fact, $C$
does not extend to $X_2$. Consequently, it is not an intersection of
the form $C=X\cap{S}$ for any surface $S\subset \bbP^4$.
\end{prop}

At this point, what would seem to be missing in a more complete
understanding of the conjecture of Griffiths and Harris is firstly,
whether the existence of the ``special'' curves of Voisin which
disprove it are indeed as special as they seem. Secondly, one would
like to know that in spite of this conjecture being false, whether
there is a ``weaker'' {\it generalised Noether-Lefschetz theorem}.

As explained below, {\it Arithmetically Cohen-Macaulay} (ACM) vector
bundles and subvarieties on hypersurfaces provide answers to both
these questions. Let $(X,\Oh_X(1))$ be a smooth polarised variety and
$\mathcal{F}$ be any coherent sheaf. Let
$\HH^i_\ast(X,\mathcal{F}):=\oplus_{\nu\in\bbZ}\HH^i(X,\mathcal{F}(\nu))$.
Recall that a vector bundle $E$ on $X$ is said to be ACM if
$\HH^i_\ast(X,E)\comment{:=\oplus_{\nu\in\bbZ}\HH^i(X,E(\nu))}=0$ for
$0< i<\dim{X}$. A subvariety $Z\subset X$ is said to be ACM if
$\HH^i_\ast(X,I_{Z/X})\comment{:=\oplus_{\nu\in\bbZ}\HH^i(X,I_{Z/X}(\nu))}=0$
for $0<i\leq \dim{Z}$. Furthermore, a codimension two subvariety
$Z\subset X$ is said to be {\it arithmetically Gorenstein}, if it is
the zero locus of a section of a rank two ACM bundle $E$ on $X$.

Let $X\subset \bbP^n$ be a smooth hypersurface. Given a codimension
two ACM subvariety $Z\subset X$ with dualising sheaf $\omega_{Z}$, one
 can associate an ACM vector bundle $E$ of rank $r+1$ where $r$ is the  minimal number 
 of generators of the $\HH^{0}_{\ast}(\bbP^{n}, \Oh_{\bbP^{n}})$-module 
 $\HH^{0}_{\ast}(Z, \omega_{Z})$.  The isomorphism 
$ \omega_{Z}\isom \sExt^{1}_{\Oh_{X}}(I_{W/X}, K_{X}),$  
where $K_{X}$ is the canonical bundle of $X$, gives rise to an isomorphism
$$ \HH^{0}(X,\oplus_{i=1}^{r}\omega_{Z}(a_{i})) \isom 
\HH^{0}(X, \sExt^{1}_{\Oh_{X}}(I_{W/X},\oplus_{i=1}^{r} K_{X}(a_{i})) )\isom 
Ext^{1}(I_{W/X}, \oplus_{i=1}^{r}K_{X}(a_{i})).$$ 
This isomorphism  takes a minimal set of generators to a rank $r+1$ ACM bundle (the fact that it 
is a bundle follows from the Auslander-Buchsbaum formula since $Y$ is locally Cohen-Macaulay)
and hence in the case  when $Z$ is arithmetically Gorenstein, this is just Serre's construction.
Conversely (see \cite{Kleiman}), given any ACM bundle $E$ of rank $r+1$ on $X$
and $r$ general sections in sufficiently high degree, one can obtain
an ACM subvariety $Z\subset X$ such that $E$ is the ACM bundle
associated to it.  In \cite{MPRV},  using the above correspondence, the following was proved:

\begin{prop}
Let $X\subset \bbP^n$ be a smooth hypersurface of degree $d>1$. If a
codimension two ACM subvariety $Z\subset X$  extends to $X_2$, then
the associated ACM vector bundle $E$ splits into a sum of line
bundles.
\end{prop}

Examples  of non-split ACM bundles on smooth hypersurfaces of
degree $d>1$ can be found  in \cite{BGS} (see \cite{MPRV} for
another construction). The existence of such bundles, together with the above proposition,
immediately implies that there exist plenty of curves in $X$, which
disprove the conjecture of Griffiths and Harris.  It can be easily
checked that Voisin's curves are in fact ACM, thus providing a
conceptual explanation why the conjecture is false.

The following theorem which can be viewed as a weak generalisation of
Theorem \ref{nlt} was proved in \cite{MPR2} and \cite{R}.

\begin{thm}\label{weaknlt}
Let $X\subset\bbP^4$ be a smooth, general hypersurface of degree
$d\geq 6$ with defining polynomial $f\in \HH^0(\bbP^4,
\Oh_{\bbP^4}(d))$. The following equivalent statements hold true.

\begin{enumerate}
\item\label{weaknlt1} Any rank two ACM bundle on $X$ is a sum of line
bundles.

\item\label{weaknlt2} $f$ cannot be expressed as the {\it Pfaffian} of
a minimal skew-symmetric matrix of size $2k\times 2k$, $2\leq k \leq
d$, whose entries are homogeneous polynomials in five variables.

\item\label{weaknlt3} A curve $C\subset X$ is a complete intersection
  if and only if $C$ is {\it arithmetically Gorenstein} i.e., $C$ is
  the zero locus of a non-zero section of a rank two ACM bundle on
  $X$.

\end{enumerate}
\end{thm}

ACM bundles on hypersurfaces have been studied earlier. To add some
history, Kleppe showed in \cite{Kleppe} that any rank two ACM bundle
$E$ on a smooth hypersurface $X\subset\bbP^n$, $n\geq 6$, splits as a
sum of line bundles.  When $n=5$, $3\leq d \leq 6$ or $n=4$ and $d=6$,
and $X$ is a general smooth hypersurface, the above splitting result
was first obtained by Chiantini and Madonna (see \cite{CM1, CM2}).
\comment{while when $n=4$ and $d=5$, the infinitesimal rigidity was first
proved in \cite{CM3}.} The first general results on ACM bundles, which
subsumed these results \comment{of \cite{Kleppe} and \cite{CM1},} were 
first proved in \cite{MPR1}. These in turn led to the proof of Theorem \ref{weaknlt},
as given in \cite{MPR2}. 
\comment{
One of the main results of that paper is the
following which we record here for future reference.

\begin{thm}
Let $X\subset\bbP^{5}$ be a smooth, general hypersurface of degree
$d\geq 3$. Then any ACM vector bundle of rank two on $X$ splits as a
sum of line bundles.
\end{thm}
}
An important ingredient of the proof of Theorem \ref{weaknlt} used in
\cite{R}, is the following theorem of Green (see \cite{MG}) and Voisin
(unpublished). 

\begin{thm}\label{aj}
Let $X\subset \bbP^4$ be a smooth, general hypersurface of degree
$d\geq 6$. Then the image of the Abel-Jacobi map 
$$\CH^2(X)_\bbQ \to J^2(X)_\bbQ$$ from the (rational) Chow group of
codimension two cycles on $X$ to the intermediate Jacobian modulo
torsion, is zero.
\end{thm}

\comment{
As mentioned earlier, Conjecture \ref{ghconj} is one in a series of
conjectures made by Griffiths and Harris. The above theorem partially
answers another in this series.
}
Notice that Noether-Lefschetz type questions can be asked more
generally for complete intersection subvarieties in projective
space. Theorem \ref{nlt} for instance, is well understood in a more
general situation (see \cite{SGA}).

\begin{thm}
Let $Y$ be a smooth projective threefold and $L$ a sufficiently ample,
base point free line bundle on $Y$. Let $X\in\vert L\vert$ be a
smooth, very general member of the linear system $\vert L\vert$. Then
the restriction map $\iota^\ast:\Pic(Y)\to \Pic(X)$ is an isomorphism.
\end{thm}

In the above theorem, we need $L$ sufficiently ample to imply that the
map $\HH^2(\Oh_Y) \to \HH^2(\Oh_X)$ is not surjective. In case
$Y=\bbP^3$, this translates to the condition $L=\Oh_{\bbP^3}(d)$ with
$d\geq 4$. In particular, this gives an extension of the
Noether-Lefschetz theorem to complete intersection surfaces of
multi-degree $(d_1,\cdots,d_{n-2})$ in $\bbP^n$ (here the corresponding
condition on the multi-degree of $X$ is $\sum_{i=1}^{n-2}d_i\geq n+1$).
Using (the infinitesimal version of) this theorem and some explicit
analysis, \comment{(in the case where $X$ is the intersection of two quadrics in
$\bbP^4$ which is the only case where $X$ is not a hypersurface and
the multi-degree does not satisfy the above inequality),} Harris and
Hulek (see \cite{HH}) extended Theorem \ref{nbs1} to smooth complete 
intersection surfaces.

\comment{
\begin{thm}
Let $X\subset \bbP^n$ be a smooth complete intersection surface of
multidegree $(d_1,\cdots,d_{n-2})$ and $C\subset X$ be a smooth
curve. Then the normal bundle sequence
$$0 \to N_{C/X} \to N_{C/\bbP^n} \to \bigoplus_{i=1}^{n-2}\Oh_C(d_i)
\to 0$$ splits if and only if $C\subset X$ is a complete intersection.
\end{thm}
}
Finally, Green and M\"uller-Stach (see \cite{G-MS}) have proved a
generalisation of Theorem \ref{aj} to complete intersection
subvarieties of sufficiently high multidegree.

\begin{thm}\label{gms}
Let $X$ be a general complete intersection subvariety in $\bbP^n$ of
sufficiently high multi-degree and dimension at least three.\comment{Let
$\CH^2(X)_\bbQ$ be the (rational) Chow group of codimension $2$
cycles modulo rational equivalence.} The image of the cycle class map
$\CH^2(X)_\bbQ \to \HH^{4}_{\mathcal{D}}(X,\bbQ)$ into Deligne
cohomology is just the image of the hyperplane class in $\bbP^n$. In
particular, the image of the Abel-Jacobi map of $X$ is contained in
the torsion points of $\JJ^2(X)$.
\end{thm}

In view of these above two theorems, it is natural to seek extensions
of Theorem \ref{weaknlt} to complete intersection subvarieties of
projective space. 

\comment{Let $X\subset\bbP^n$ be a complete intersection with multi-degree
$(d_1,\cdots,d_k)$. If $d_i\gg 0$ $\forall$ $1\leq i \leq k$, then we
say that $X$ is a complete intersection of {\it sufficiently large
  multi-degree}. 
}

\section{Main results}
 The main results of this note are the following.

\begin{thma}\label{hdims}
Let $X\subset \bbP^n$ be a general smooth complete intersection
subvariety of dimension at least four and sufficiently high
multidegree. Then any ACM vector bundle of rank two on $X$ is a direct
sum of line bundles.
\end{thma}

This can be viewed as a generalisation to complete intersections of
the main result of \cite{MPR1}.

Let $d_1\leq d_2\leq \cdots\leq d_{n-3}$ be a sequence of positive
integers, with $d_i\gg 0$ for $1\leq i \leq n-3$ and $d_{n-3}>
\mbox{max}\{d_i+d_j\}$ for $1\leq i,~j \leq n-4$. Recall that a bundle
$E$ is said to be {\it normalised} if $h^0(E(-1))=0$ but $h^0(E)\neq
0$. With this notation, we have

\begin{thmb}\label{3fold}
Let $X\subset \bbP^n$ be a general smooth complete intersection
threefold of sufficiently high multidegree $(d_1,\cdots,d_{n-3})$ as
above. Then any arithmetically Cohen-Macaulay, normalised rank two
vector bundle $E$ on $X$ is a direct sum of line bundles provided
$c_1(E)< d_{n-3}-1$.
\end{thmb}

Theorem A is obtained as a consequence of Theorem B. A word about the
inequality satisfied by the first Chern class and about the condition
$d_i\gg 0$ in the above theorems: since any bundle $E$ splits iff
$E(m):=E\tensor\Oh_X(m)$ splits for some $m$, we may assume that $E$
is normalised. Madonna showed in \cite{Madonna} that on any smooth
three dimensional complete intersection $X\subset \bbP^n$ a normalised
rank two ACM bundle splits as a direct sum of two line bundles unless
$-\sum_{i=1}^{n-3}d_i+n-1\leq c_1(E) \leq \sum_{i=1}^{n-3}d_i+3-n$. We
notice that when $n=4$ and $X$ is a general hypersurface of degree
$d\gg0$, combining our result and Madonna's bound, it follows that the
only possible first Chern class of a normalised and indecomposable
rank two ACM bundle on $X$ is $d-1$. This case is classical. Indeed,
the existence of an indecomposable ACM bundle of rank two on a general
smooth hypersurface in $\bbP^4$ is equivalent to the fact that a
general homogeneous polynomial in five variables can be obtained as
the Pfaffian of a (minimal) skew-symmetric matrix of linear forms. It
is well known (see \cite{Beauville}), by a simple dimension count,
that when $d\geq 6$, this is not possible. Thus Theorem B can be seen
as a generalisation of Theorem \ref{weaknlt} above. Finally, to prove
Theorem B, we argue as in \cite{R} for the cases of general three
dimensional hypersurfaces of degree at least six. Indeed we use here
the result of Green and M\"uller-Stach (Theorem \ref{gms}) which
generalised the corresponding result of Green and Voisin (Theorem
\ref{aj}). In doing so, we need to restrict to the cases of complete
intersections of sufficiently high multi-degree and dimension at least $3$.
 
We have then the following interesting:

\begin{corc}
Let $X\subset \bbP^n$ be a general smooth complete intersection
subvariety of sufficiently high multidegree.
\begin{enumerate}
\item Any arithmetically Gorenstein subvariety $T\subset X$ of
  codimension two is a complete intersection in $X$ provided
  $\dim{X}\geq 4$.
\item Suppose $X$ is a threefold, and $C\subset X$ is any
  arithmetically Gorenstein curve. Then $C$ is the intersection of $X$
  with a codimension two subscheme $S\subset \bbP^n$ if and only if
  $C$ is a complete intersection in $X$. In addition, if the
  rank two bundle $E$ associated to $C$ via Serre's correspondence is normalised,
  and $c_1(E)<d_{n-3}-1$, then $C$ is a complete intersection in $X$.
\end{enumerate}
\end{corc}

Finally, we should mention that another motivation for the questions
on ACM vector bundles comes from the conjectures of
Buchweitz-Gruel-Schreyer (see \cite{BGS}, Conjecture B) on the
triviality of low rank ACM bundles on hypersurfaces (see \cite{R} for
more details). 

The present paper builds on results proved and techniques developed in
\cite{MPR1, MPR2, R, Wu1, Wu2}, some of which have been included here
for the sake of completeness. The non-degeneracy of the infinitesimal
invariant in particular, is shown by refining Xian Wu's proof in
\cite{Wu1}.

\section{Acknowledgements}
The authors would like to thank the referee for comments which helped make the paper more readable.

\section{The infinitesimal invariant associated to a normal 
function}\label{iinf}

Let $Y$ be a smooth projective variety of dimension $2m$, $m\geq 1$.
Let $\mathcal{X} \to S$ be the universal family of smooth, degree $d$
hypersurfaces of in $Y$. Let $\mathcal{C}\subset\mathcal{X}$ be a
family of codimension $m$ subvarieties over $S$. If $l=l(s)$
is the degree of $C_s$ and $D_s$ is a codimension $m$ linear section
(i.e. $D_s$ is an intersection of $m$ hyperplanes in $X_s$), then the
family of cycles $\mathcal{Z}$ with fibre
$\mathcal{Z}_s:=d~\mathcal{C}_s-lD_s$ for $s\in S$ defines a
fibre-wise null-homologous cycle, i.e. an element in
$\CH^m(\mathcal{X}/S)_{hom}$. Let $\mathcal{J}:=\{J(X_s)\}_{s\in S}$
be the family of intermediate Jacobians. In such a situation,
Griffiths (see \cite{G}) has defined a holomorphic function
$\nu_{\mathcal{Z}}:S \to \mathcal{J}$, called the {\it normal
  function}, which is given by $\nu_{\mathcal{Z}}(s)=\mu_s(Z_s)$ where
$\mu_s:\CH^m(X_s)_{\mbox{\small{hom}}} \to J(X_s)$ is the {\it
  Abel-Jacobi} map from the group of null-homologous cycles to the
intermediate Jacobian.

This normal function satisfies a ``quasi-horizontal'' condition (see
\cite{Vbook}, Definition 7.4).  Associated to the normal function
$\nu_{\mathcal{Z}}$ above, Griffiths (see \cite{G1} or \cite{Vbook}
Definition 7.8) has defined the infinitesimal invariant
$\delta\nu_{\mathcal{Z}}$. Later Green \cite{MG} generalised this
definition and showed that Griffiths' original infinitesimal invariant
is just one of the many infinitesimal invariants that one can
associate to a normal function.  For a point $s_0\in S$, let
$X=X_{s_0}$, $C:=\mathcal{C}_{s_0} \subset X$ and $D:=D_{s_0}$. Green
showed that in particular $\delta\nu_{\mathcal{Z}}(s_0)$ is an element
of the dual of the middle cohomology of the following (Koszul) complex
\begin{equation}\label{koszul}
\ext{2}\HH^1(X,T_X)\tensor\HH^{m+1,m-2}(X)\to
  \HH^1(X,T_X)\tensor\HH^{m,m-1}(X)\to
  \HH^{m-1,m}(X).
\end{equation}
We now specialise to the case $m=2$ where $X\subset Y$ is a smooth
hypersurface of dimension $3$ and $C\subset X$ is a curve of degree
$l$. Then $Z:=dC-lD$ is a nullhomologous $1$-cycle with support
$W:=C\bigcup D$.  
At a point $s\in S$, this infinitesimal invariant is therefore a functional
$$\delta\nu_{\mathcal{Z}}(s):
\ker\left(\HH^1(X,T_X)\tensor\HH^1(X,\Omega^2_X) \to
\HH^2(X,\Omega^1_X)\right)\to\bbC.$$ 
Consider the composite map 
$$\gamma: \HH^1(X,T_X)\tensor\HH^1(X,\mathcal{I}_{W/X}\tensor\Omega^2_X) \to 
\HH^1(X,T_X)\tensor\HH^1(X,\Omega^2_X) \to \HH^2(X,\Omega^1_X).$$
By abusing notation, we will let 
$$\delta\nu_{\mathcal{Z}}(s):\ker\gamma\to \bbC$$ 
denote the composite map.
On the other hand, starting with the short exact sequence
$$ 0 \to  I_{W/X}\to \Oh_{X} \to \Oh_{W} \to 0,$$
and tensoring with $\Omega^{1}_{X}$, yields a
long exact sequence of cohomology
$$ \cdots \to \HH^1(\Omega^1_X) \to \HH^1(\Omega^1_X\tensor\Oh_W)
\to\HH^2(I_{W/X}\tensor\Omega^1_X) \to 
\HH^2(\Omega^1_X) 
\to  0 .$$ 
Combining this sequence with the Koszul complex \eqref{koszul},  
we get a commutative diagram:
\begin{equation}\label{formula}
\begin{diagram}{ccccccccc}
& & & & \HH^1( T_X)\tensor\HH^1(I_{W/X}\tensor\Omega_X^2) & &
& &
\\
& & & & \downarrow{\beta} & \searrow{\gamma}& & & \\
0 & \to & \frac{\HH^1(\Omega^1_X\tensor\Oh_W)}{\HH^1(\Omega^1_X)}
&
\by{\lambda} &\HH^2(I_{W/X}\tensor\Omega^1_X) & \to &
\HH^2(\Omega^1_X) &
\to & 0 \\
& & \downarrow{\chi} & & & &  & & \\
& & \bbC & & & & & & \\
\end{diagram}
\end{equation}
where $\chi$ is given by integration over the cycle $Z$ and $\beta$ is
a Koszul map. As a result, one has an induced map
$$\ker\gamma \to  \frac{\HH^1(\Omega^1_X\tensor\Oh_W)}{\HH^1(\Omega^1_X)} .$$
The following is the main result that we shall use in this paper.
\begin{thm}[Griffiths \cite{G1,G2}]
Let $\nu_{\mathcal{Z}}$ be the normal function as described
above. \comment{Consider the following diagram:
\begin{equation}
\begin{diagram}{ccccccccc}
& & & & \HH^1( T_X)\tensor\HH^1(I_{W/X}\tensor\Omega_X^2) & &
& &
\\
& & & & \downarrow{\beta} & \searrow{\gamma}& & & \\
0 & \to & \frac{\HH^1(\Omega^1_X\tensor\Oh_W)}{\HH^1(\Omega^1_X)}
&
\by{\lambda} &\HH^2(I_{W/X}\tensor\Omega^1_X) & \to &
\HH^2(\Omega^1_X) &
\to & 0 \\
& & \downarrow{\chi} & & & &  & & \\
& & \bbC & & & & & & \\
\end{diagram}
\end{equation}
where $\chi$ is given by integration over the cycle $Z$.} Then 
$\delta\nu_{\mathcal{Z}}(s_0)$, the infinitesimal invariant evaluated at
the point $s_0\in S$, is the composite
$$ \ker{\gamma}\to \frac{\HH^1(\Omega^1_X\tensor\Oh_W)}
{\HH^1(\Omega^1_X)}\by{\chi} \bbC.$$
\end{thm}

\comment{We identify the tangent space $T_{s}=\HH^0(X,\Oh_X(d))$ and
$\HH^1(I_{W/X}\tensor T_X)$ with their respective images in $\HH^1(
T_X)$.} 
\begin{remark}
The map $\chi$ can be understood as follows: Since $D$ is a
general plane section of $X$, by Bertini $C\cap{D}=\emptyset$.  Thus
$\Oh_W\isom\Oh_C\oplus\Oh_D$ and so
$$\HH^1(\Omega^1_X\tensor\Oh_W)\isom \HH^1(\Omega^1_X\tensor\Oh_C)
\oplus \HH^1(\Omega^1_X\tensor\Oh_D). $$ For any irreducible curve
$T\subset X$, let $r_T:\HH^1(\Omega^1_X\tensor\Oh_T) \to
\HH^1(\Omega^1_T)\isom \bbC$ be the natural restriction map.
 For any element $(a,b)\in
\HH^1(\Omega^1_X\tensor\Oh_W)$, $\chi(a,b):=dr_C(a)-lr_D(b) \in
\bbC$. Clearly, this map factors via the quotient
$\HH^1(\Omega^1_X\tensor\Oh_W)/\HH^1(\Omega^1_X).$
\end{remark}

The main result in this situation is Theorem \ref{gms} which implies in 
particular that the normal function is zero on an open subset of the 
parameter space. By Theorem 1.1 in \cite{MG}  (or Proposition 1.2.3 of \cite{Wu1}), 
we then have the following.

\begin{thm}\label{gmseq}
Let $X$ be a general complete intersection subvariety in $\bbP^n$ of
sufficiently high multi-degree and dimension at least
three.
\comment{Let $\CH^2(X)\tensor\bbQ$ be the (rational) Chow group
  of codimension $2$ cycles modulo rational equivalence. The image of
  the cycle class map $\CH^2(X)\tensor\bbQ \to
  \HH^{4}_{\mathcal{D}}(X,\bbQ)$ into the (rational) Deligne cohomology
  is just the image of the hyperplane class in $\bbP^n$. In
  particular, the image of the Abel-Jacobi map of $X$ is contained in
  the torsion points of $\JJ^2(X)$. Equivalently,} If $\mathcal{Z}\to
S$ is a family of codimension two, degree zero cycles contained in the
universal complete intersection $\mathcal{X}\subset \bbP^n\times S$,
then the infinitesimal invariant $\delta\nu_{\mathcal{Z}}$ associated
to $\mathcal{Z}$ vanishes at a general point $s\in S$.
\end{thm}

\section{ACM bundles on a smooth subvariety $X\subset Y$}
Let $X=\bigcap_{i=1}^{n-3}Y_i$ be a general complete intersection of
smooth hypersurfaces $Y_i\subset\bbP^n$ of degree $d_i$.  Let $E\to X$
be an indecomposable, normalised ACM bundle of rank two. In this
section, we shall establish several lemmas which will enable us to
prove the non-degeneracy of the infinitesimal invariant coming from a
family of arithmetically Gorenstein curves. The criterion is a
refinement of Wu's criterion (see \cite{Wu1}).

\begin{lemma}\label{c1inequality}
Let $E$ be as above with first Chern class $\alpha$. Then the zero
locus of every non-zero section of $E$ has codimension $2$ in $X$. If
$C\subset X$ is the zero locus of a section of $E$, then we have the
exact sequence
\begin{eqnarray}\label{serre}
0\to \Oh_X(-\alpha) \to E^\vee \to I_{C/X} \to 0.
\end{eqnarray}
Furthermore, $E$ is $((\sum_{i=1}^{n-3}d_i)-n+3-\alpha)$--regular,
$n-1-\sum_{i=1}^{n-3}d_i\leq \alpha\leq (\sum_{i=1}^{n-3}d_i)-n+3$, and
$K_C=\Oh_C(\sum_{i=1}^{n-3}d_i-n-1+\alpha)$.
\end{lemma}

\begin{proof}
See \cite{R} for the proof of the first part of the Lemma. The
regularity $\rho$ of $E$ can be computed easily (see {\it
  op.~cit.}). For the inequality satisfied by $\alpha$, see
\cite{Madonna}.  \comment{ Thus the upper bound for $\rho$ is easily
  seen. For the lower bound and the rest of the proof, see proof of
  Lemma 3.3 in \cite{MPR1}.}
\end{proof}

\noindent{\bf Notation:} Let $d_1\leq d_2\leq \cdots\leq d_{n-3}$ be a 
sequence of sufficiently large positive integers (i.e. $d_i\gg 0$), and 
$d_{n-3}>\mbox{max}\{d_i+d_j\}$ for $1\leq i,~j\leq n-4$. For the rest of the
paper, we will denote by $$Y:=\bigcap_{i=1}^{n-4}Y_i, \hspace{2mm}
d=d_{n-3} \hspace{3mm} \mbox{and} \hspace{3mm} X\in|\Oh_Y(d)|$$ a
general member. Also, we will let $E$ denote a normalised, indecomposable ACM
bundle of rank two on $X$ and $C\subset X$ to denote the zero locus of a
non-zero section of $E$ (which is a curve by Lemma
\ref{c1inequality}).
\subsection{Towards the surjectivity of the map $\chi$}
The main result of this subsection is the surjectivity of the map $\chi$.  
This is achieved by identifying a subspace of $\HH^{1}(W, \Omega^{1}_{X|W})$,
restricted to which $\chi$ is surjective (see Corollary \ref{chiisonto}). The proof 
crucially depends on the following:
\begin{lemma}\label{firststep}
Let $C\subset X\subset Y\subset\bbP^{n}$ be as above. Then the natural map
$$\HH^{1}(X,\Omega^{2}_{\bbP^{n}}(d)_{|X}) \to\HH^{1}(C,\Omega^{2}_{\bbP^{n}}(d)_{|C}) $$
is zero.
\end{lemma}

Before we prove this lemma, we shall need several results which we shall prove now. Let
\begin{equation}\label{ses1}
0 \to G_Y \to F_{0,Y}\to E \to 0, 
\end{equation}
be a minimal resolution of $E$ by vector bundles on $Y$. So
$F_{0,Y}:=\bigoplus\Oh_Y(-a_i)$, where $a_i\geq 0$ and the kernel
$G_Y$ is ACM. The fact that $G_Y$ is a bundle follows from the
Auslander-Buchsbaum formula (see \cite{Eisenbud}, Chapter
19). 

Applying $\sHom_{\Oh_Y}(\cdot\,,\Oh_Y)$ to sequence (\ref{ses1}),
we get 
\begin{equation}\label{dualres}
0\to F_{0,Y}^\vee \to G_Y^\vee \to E^\vee(d) \to 0.
\end{equation}
(see \cite{MPR1} where it is proved when $X\subset\bbP^{n}$ is a hypersurface: 
the same argument works on replacing $\bbP^{n}$ by $Y$)

\begin{lemma}[see also \cite{R}]\label{split}
There exists a commutative diagram
\begin{equation}
\begin{array}{ccccccccc} \label{c1}
0 & \to & F_{0,Y}^\vee & \to & G_Y^\vee & \to & E^\vee(d)
& \to & 0 \\
& & \downarrow{\phi}  &  & \downarrow & & \downarrow{s^\vee} &  &  \\
0 & \to & \Oh_{Y} & \to & \Oh_{Y}(d) & \to & \Oh_{X}(d) & \to & 0
\\
\end{array}
\end{equation}
where under the isomorphism $$\Hom(F_{0,Y}^\vee,
\Oh_Y)\isom\HH^0(F_{0,Y})\isom\HH^0(E)
\comment{\isom\Ext^1_{\Oh_Y}(E^\vee(d),\Oh_Y)}, \hspace{3mm} \phi
\mapsto s. $$ In addition, $F_{0,Y}^\vee\by{\phi}\Oh_Y$ is a split
surjection (i.e. $\phi$ is the projection onto one of the factors).
\end{lemma}

\begin{proof}
We consider the following push-out diagram:
\begin{equation}
\begin{array}{ccccccccc} 
0 & \to & F_{0,Y}^\vee & \to & G_Y^\vee & \to & E^\vee(d)
& \to & 0 \\
& & \downarrow{\phi}  &  & \downarrow & & || &  &  \\
0 & \to & \Oh_{Y} & \to & \mathcal{K} & \to & E^\vee(d) & \to & 0
\\
\end{array}
\end{equation}

Since $F_{0,Y}^\vee=\oplus\Oh_Y(a_i)$, $a_i\geq 0$, any such diagram
corresponds to a section $\phi\in\HH^0(F_{0,Y})$.

Next consider the pull-back diagram:
\begin{equation}
\begin{array}{ccccccccc} 
0 & \to & \Oh_{Y} & \to & \mathcal{K} & \to & E^\vee(d) & \to & 0\\
\comment{0 & \to & F_{0,Y}^\vee & \to & G_Y^\vee & \to & E^\vee(d)
& \to & 0 \\}
& & ||  &  & \downarrow & & \downarrow{s^\vee} &  &  \\
0 & \to & \Oh_{Y} & \to & \Oh_{Y}(d) & \to & \Oh_{X}(d) & \to & 0
\\
\end{array}
\end{equation}
Any such diagram corresponds to a section $s\in\HH^0(E)$. Since
$\HH^0(F_{0,Y})\isom \HH^0(E)$, there is a bijective correspondence
between the diagrams above. Combining them, we get the desired
commutative diagram. The morphism $\phi$ is a split surjection since
$a_i\geq 0,$ $\forall~ i$.
\end{proof}

Tensoring the exact sequence (\ref{dualres}) with $\Omega^2_{\bbP^n}$
and taking cohomology, we get 
\begin{equation}\label{cohdualres}
\to  \HH^1(E^\vee(d)\tensor\Omega^2_{\bbP^n})  \to
\HH^2(F^\vee_{0,Y}\tensor\Omega^2_{\bbP^n})  \to 
\HH^2(G^\vee_Y\tensor\Omega^2_{\bbP^n}) \to 
\end{equation}

\begin{lemma}\label{mapiszero}
The map $\HH^2(F^\vee_{0,Y}\tensor\Omega^2_{\bbP^n}) \to
\HH^2(G^\vee_Y\tensor\Omega^2_{\bbP^n})$ in diagram (\ref{cohdualres}) is the
zero map.
\end{lemma}

\begin{proof}
Let $F_1\to G_Y^\vee$ be a surjection from a sum of line bundles on
$\bbP^n$ to $G_Y^\vee$, induced by a minimal set of generators of
$G_Y^\vee$. Let $F_0$ be a sum of line bundles on $\bbP^n$ such that
$F_0\tensor\Oh_Y=F_{0,Y}$. \comment{Then
$$\HH^2(F^\vee_{0}\tensor\Omega^2_{\bbP^n}) \isom
  \HH^2(F^\vee_{0,Y}\tensor\Omega^2_{\bbP^n}).$$} The map
$F_{0,Y}^\vee \to G_Y^\vee$ lifts to a map $\Phi:F_0^\vee \to F_1$,
since $G_Y^\vee$ is ACM. Hence we have a commuting square

\[
\begin{array}{ccc} \label{c3}
\HH^2(F^\vee_{0}\tensor\Omega^2_{\bbP^n}) & \to &
\HH^2(F_{1}\tensor\Omega^2_{\bbP^n}) \\
 \downarrow{\isom}  &  & \downarrow \\
\HH^2(F^\vee_{0,Y}\tensor\Omega^2_{\bbP^n}) & \to &
\HH^2(G^\vee_Y\tensor\Omega^2_{\bbP^n}) \\
\end{array}
\]

To prove the lemma, it is enough to prove that the top horizontal map,
which is given by the matrix $\Phi$, is zero. In other words, we need
to show that $\Phi$ has no non-zero scalar entries.

Suppose there was such an entry, then we would have a diagram

\[
\begin{array}{ccccccccc} \label{conmin}
 &  & 0 &  & 0 &  & 
&  &  \\
& & \downarrow  &  & \downarrow & &  &  &  \\
 &  & \Oh_Y & = & \Oh_Y &  & 
&  &  \\
& & \downarrow  &  & \downarrow & &  &  &  \\
0 & \to & F_{0,Y}^\vee & \to & G_Y^\vee & \to & E^\vee(d)
& \to & 0 \\
& & \uparrow\downarrow  &  & \downarrow & & || &  &  \\
0 & \to & \bar{F} & \to & \bar{G} & \to & E^\vee(d)
& \to & 0 \\
& & \downarrow  &  & \downarrow & &  &  &  \\
 &  & 0 &  & 0 &  & 
&  &  \\
\end{array}
\]

Here $\bar{F}$ (resp. $\bar{G}$) is defined as the cokernel of the
inclusion $\Oh_Y\into F_{0,Y}^\vee$ (resp. $\Oh_Y\into G_{Y}^\vee$).
Applying $\sHom_{\Oh_Y}(.\, ,\Oh_Y)$ to the diagram above, we get

\[
\begin{array}{ccccccccc} \label{dualconmin}
 &  & 0 &  & 0 &  & 
&  &  \\
& & \downarrow  &  & \downarrow & &  &  &  \\
0 & \to & \bar{G}^\vee & \to & \bar{F}^\vee & \to & 
Ext^1_{\Oh_Y}(E^\vee(d),\Oh_Y)=E
& \to &  \\
& & \downarrow  &  & \downarrow & & || &  &  \\
0 & \to & G_Y & \to & F_{0,Y} & \to & 
Ext^1_{\Oh_Y}(E^\vee(d),\Oh_Y)=E
& \to & 0 \\
& & \downarrow  &  & \uparrow\downarrow & &  &  &  \\
 &  & \Oh_Y & = & \Oh_Y &  & 
&  &  \\
& & \downarrow  &  & \downarrow & &  &  &  \\

 &  &  &  & 0 &  & 
&  &  \\
\end{array}
\]

Since $\Oh_Y$ is a summand of $F_{0,Y}$ which is in the image of the
map $G_Y \to F_{0,Y}$, the composite map $\Oh_Y \to F_{0,Y} \to E$ is
zero. In particular, this implies that $G_Y\to \Oh_Y$ is a surjection
and so
$$0 \to \bar{G}^\vee \to \bar{F}^\vee \to E \to 0 $$ is also a
resolution. This contradicts the minimality of sequence (\ref{ses1}).
\end{proof}

\begin{proof}[Proof of Lemma \ref{firststep}]
We have the following resolution for $\Oh_X$ on $\bbP^n$:
$$0 \to \Oh_{\bbP^n}(-\sum_{i=1}^{n-3}d_i) \to \cdots \to
\bigoplus_{i=1}^{n-3}\Oh_{\bbP^n}(-d_i) \to \Oh_{\bbP^n}\to \Oh_{X}
\to 0.$$ Using this resolution and the fact that
$d=d_{n-3}>\mbox{max}_{1\leq i,j\leq n-4}\{d_i+d_j\}$, one can show that
$$\HH^1(\Oh_{X}(d)\tensor\Omega^2_{\bbP^n})\isom \bigoplus_{i=1}^{n-3}
\HH^2(\Omega^2_{\bbP^n}(d-d_i))\isom\HH^2(\Omega^2_{\bbP^n}).$$
Tensoring diagram (\ref{c1}) with $\Omega^2_{\bbP^n}$ yields a
commuting square
\begin{equation}
\begin{array}{ccccccccc} \label{c2}
& & \HH^1(E^\vee(d)\tensor\Omega^2_{\bbP^n}) & \onto &
\HH^2(F^\vee_{0,Y}\tensor\Omega^2_{\bbP^n}) & & & & \\
& & \downarrow & & \downarrow & & & & \\
 & &
\HH^1(\Oh_{X}(d)\tensor\Omega^2_{\bbP^n}) & \isom &
\HH^2(\Oh_{Y}\tensor\Omega^2_{\bbP^n}) &  & & & \\
\end{array}
\end{equation}

Here the top horizontal map is a surjection by Lemma \ref{mapiszero},
and the right vertical arrow is a surjection by Lemma \ref{split}.
Hence the map
$\HH^1(E^\vee(d)\tensor\Omega^2_{\bbP^n})\to
\HH^1(\Oh_{X}(d)\tensor\Omega^2_{\bbP^n})$ is a surjection. Since this
map factors via $\HH^1(I_{C/X}(d)\tensor\Omega^2_{\bbP^n})$, the
natural map $\HH^1(I_{C/X}(d)\tensor\Omega^2_{\bbP^n})\to
\HH^1(\Oh_{X}(d)\tensor\Omega^2_{\bbP^n})$ is a surjection. Hence we
are done.
\end{proof}

Let $h_Y\in\HH^1(\Omega^1_Y)$ be the restriction of the generator
$h\in\HH^1(\Omega^1_{\bbP^n})$ and consider the class
$h_Y^2\in\HH^2(\Omega^2_Y)$. This is the image of $h_X$, the
hyperplane class in $X$, under the Gysin map $\bbC\isom
\HH^1(\Omega^1_X) \to\HH^2(\Omega^2_Y)$.  Furthermore, since $d\gg0$,
by Serre vanishing, we have $\HH^i(\Omega^2_Y(d))=0$ for $i=1,2$. Hence the
coboundary map $\HH^1(\Omega^2_Y(d)_{|X})\to\HH^2(\Omega^2_Y)$ is an
isomorphism.  By abuse of notation, we will denote the inverse image  of $h^{2}_{Y}$ 
under this isomorphism by $h^{2}_{Y}$.
\begin{cor}
Under the natural map
$$\HH^1(\Omega^2_Y(d)_{|X}) \to \HH^1(\Omega^2_Y(d)_{|C}), \hspace{5mm}
h_Y^2\mapsto 0.$$ 
\end{cor}

\begin{proof}
One has a commutative square
\[
\begin{array}{ccccc} 
\HH^1(I_{C/X}(d)\tensor\Omega^2_{\bbP^n}) & \onto &
\HH^1(\Oh_X(d)\tensor\Omega^2_{\bbP^n}) & \to &
\HH^1(\Oh_C(d)\tensor\Omega^2_{\bbP^n}) \\ 
\downarrow & & \downarrow & & \downarrow \\
\HH^1(I_{C/X}(d)\tensor\Omega^2_{Y}) & \to &
\HH^1(\Oh_X(d)\tensor\Omega^2_{Y}) & \to & 
\HH^1(\Oh_C(d)\tensor\Omega^2_{Y})\\
\end{array}
\]

The first horizontal arrow in the top row is surjection, and the  middle
vertical map $\HH^1(\Oh_X(d)\tensor\Omega^2_{\bbP^n})\to
\HH^1(\Oh_X(d)\tensor\Omega^2_{Y})$ can be identified with the map
$\HH^2(\Omega^2_{\bbP^n}) \to \HH^2(\Omega^2_Y)$ which takes the
element $h^2\mapsto h_Y^2$. Hence $h_Y^2 \mapsto 0$ under the map
$\HH^2(\Omega^2_Y)\isom\HH^1(\Omega^2_Y(d)_{|X}) \to
\HH^1(\Omega^2_Y(d)_{|C})$.
\end{proof}

Now we are in  a position to prove the first step i.e., the surjectivity of the map $\chi$.
Consider the exact sequence
$$0 \to \Oh_X(-d) \to \Omega^1_{Y{|X}} \to \Omega^1_X \to 0.$$
Taking second exteriors and tensoring the resulting sequence by
$\Oh_X(d)$, we get a short exact sequence
\begin{equation}\label{2ndext}
0 \to \Omega^1_X \to \Omega^2_{Y}(d)_{|X}\to\Omega^2_X(d)\to 0.
\end{equation}

For the inclusion $C\subset X$, the natural map $\Omega^1_{X{|C}} \to
\Omega^1_C$ yields a push out diagram:
\[
\begin{diagram}{ccccccccc}
0 & \to & \Omega^1_{X{|C}} & \to & \Omega^2_{Y}(d)_{|C} &
\to &
\Omega^2_X(d)_{|C} & \to & 0 \\
& & \downarrow & & \downarrow & & || & & \\
0 & \to & \Omega^1_C & \to &
\mathcal{F} & \to & \Omega^2_X(d)_{|C} & \to & 0. \\
\end{diagram}
\]
where $\mathcal{F}$ is defined by the diagram.
\begin{lemma}
The map $\HH^1(C,\Omega^1_C) \to \HH^1(C,\mathcal{F})$ in the
associated cohomology sequence of the bottom row in the above diagram
is zero. Thus we have a surjection
$$\VV_C:=\ker[\HH^1(\Omega^1_{X|C}) \to
\HH^1(\Omega^2_{Y}(d)_{|C})] \onto \HH^1(C,\Omega^1_C).$$
\end{lemma}

\begin{proof}
We have a commutative diagram
\[\begin{array}{ccc}
\HH^1(\Omega^1_X) & \to & \HH^1(\Omega^2_{Y}(d)_{|X}) \\
\downarrow &  & \downarrow \\
\HH^1(\Omega^1_{X{|C}}) & \to &
\HH^1(\Omega^2_{Y}(d)_{|C}) \\
\downarrow &  & \downarrow \\
\HH^1(\Omega^1_C) & \to & \HH^1(\mathcal{F}) \\
\end{array}\]
The composite of the vertical maps on the left is the map which takes the class
$h_{X} \mapsto h_{C}$. Since these are the respective generators of these cohomology groups both of which are one dimensional, this composite is an isomorphism. On the other hand, the composite
$$\HH^1(\Omega^1_X) \to \HH^1(\Omega^2_Y(d)_{|X}) \to
\HH^1(\Omega^2_Y(d)_{|C})$$ is zero:  this is because the map 
$\HH^1(\Omega^1_X) \to \HH^1(\Omega^2_Y(d)_{|X})$ can be
identified with the Gysin map $\HH^1(\Omega^1_X)\to\HH^2(\Omega^2_Y)$, 
and so by the above Corollary, the generator $h_{X} \mapsto 0$ under the composite.
This implies that the map $\HH^1(C,\Omega^1_C) \to
\HH^1(C,\mathcal{F})$ is zero and so we have a surjection
$\VV_C \onto \HH^1(C,\Omega^1_C)$.
\end{proof}

\begin{cor}[Surjectivity of $\chi$]\label{chiisonto}
The composite map 
\[
\begin{array}{ccc}
\VV_C \into \ker[\HH^1(\Omega^1_{X|W}) \to
\HH^1(\Omega^2_{Y}(d)_{|W})] & \stackrel{\chi}\to & \bbC \\
\end{array}
\]
is a surjection. Hence $\chi$ is a surjection.
\end{cor}

\begin{proof}
This first inclusion follows from the fact that $\Oh_W \isom \Oh_C
\oplus \Oh_D$. The surjectivity of the composite follows from the
definition of $\chi$ and the above lemma.
\end{proof}

\subsection{Some vanishing lemmas} In this subsection, we shall prove vanishing of 
certain cohomologies. The technical condition in Theorem B is required for these
 vanishings to hold and that is the only reason for its appearance in in the statement 
 of the theorem. The main result here is Lemma \ref{tgtvanishing} and the reader 
 may skip the details which are pretty standard arguments if s/he so wishes. 
 
\begin{lemma}\label{vanishings}
With notation as above and $\alpha<d-1$, we have 
$$\HH^j(T_Y\tensor K_Y\tensor{I}_{C/X}(2d-j))=0, \hspace{3mm} j=1,~2.$$
\end{lemma}

\begin{proof}
From the exact sequence $$0 \to \Oh_X(-\alpha)\to E^\vee \to I_{C/X}
\to 0,$$ it is enough to prove the following
\begin{enumerate}
\item
$\HH^{j+1}(T_Y\tensor K_Y\tensor\Oh_X(2d-j-\alpha))=0 \hspace{3mm} j=1,~2.$
\item
$\HH^j(T_Y\tensor K_Y\tensor{E}(2d-j-\alpha))=0 \hspace{3mm} j=1,~2.$
\end{enumerate}
The vanishings in (1) above follow, on tensoring the exact sequence 
$$0 \to T_Y\tensor\Oh_X \to T_{\bbP^n}\tensor\Oh_X\to
\oplus_{i=1}^{n-4}\Oh_X(d_i) \to 0,$$ with $K_Y\tensor\Oh_Y(2d-\alpha-j)$
and using the vanishing of the following terms
\begin{itemize}
\item[(A)] $\HH^j(K_Y(d_i+2d-\alpha-j)_{|X})$ for $j=1,~2$. Since
$K_Y\isom \Oh_Y(\sum_{i=1}^{n-4}d_i-n-1)$ and $X$ is a complete
intersection, this follows.

\item[(B)] $\HH^{j+1}(T_{\bbP^n}\tensor K_Y(2d-\alpha-j)_{|X})$ for $j=1,~2$.
\end{itemize}

Using the Euler sequence restricted to $X$:
$$0\to \Oh_{X} \to \Oh_{X}(1)^{\oplus{n+1}} \to T_{\bbP^n}\tensor\Oh_X
\to 0,$$ the vanishing in (B) follows from the vanishing of
\begin{itemize}
\item $\HH^{j+1}(K_Y(2d-\alpha-j+1)_{|X})$ for $j=1,~2$ and 
\item $\HH^{j+2}(K_Y(2d-\alpha-j)_{|X})$ for $j=1,~2$.
\end{itemize}
The only non-trivial cases are when $j=2$ in the first case and $j=1$
in the second case. These vanishings hold provided $\alpha<d-1$.

For the vanishings in (2), we use the minimal resolution of $E$ on
$Y$:
$$0 \to G_Y \to F_{0,Y} \to E \to 0.$$
Then it suffices to show that

\begin{enumerate}
\item[(C)]$\HH^{j}(T_{Y}\tensor K_Y\tensor F_{0,Y}(2d-\alpha-j))=0$
for $j=1,~2$. 

\item[(D)]$\HH^{j+1}(T_{Y}\tensor K_Y\tensor
G_Y(2d-\alpha-j))$ for $j=1,~2$.
\end{enumerate}

For (C): $F_{0,Y}=\bigoplus\Oh_Y(-a_i)$ where $-a_i+\reg{E}\geq 0$. So
the above term is $\bigoplus_i \HH^{j}(T_{Y}(b_i))$ where $b_i\geq
d-j-4$. Since $d>>0$, this is true by Serre vanishing.

For (D): From the tangent bundle sequence
$$0 \to T_Y \to T_{\bbP^n}\tensor\Oh_Y \to
\oplus_{i=1}^{n-4}\Oh_Y(d_i)\to 0,$$ the required vanishings follow
from

$\bullet$ $\HH^{j}(\Oh_{Y}(d_i)\tensor K_Y\tensor G_Y(2d-\alpha-j))=0$ for 
$j=1,~2$ since $G_Y$ is ACM, and

$\bullet$ $\HH^{j+1}(T_{\bbP^n}\tensor K_Y\tensor G_Y(2d-\alpha-j))=0$
for $j=1,~2$. For this use the Euler sequence to reduce this statement
to vanishing like the above and then use the fact that $G_Y$ is ACM.
\end{proof}
\comment{
\begin{remark}
The above vanishing is the only reason why we have the technical condition 
$\alpha<d-1$ in the statement of Theorem B.
\end{remark}
}
Recall that $W=C\cup D$.
\begin{lemma}\label{tgtvanishing}
$\HH^1(T_Y\tensor K_Y\tensor\mathcal{I}_{W/Y}(2d))=0.$
\end{lemma}

\begin{proof}
From the exact sequence 
$$0\to \Oh_Y(-d) \to \mathcal{I}_{W/Y}\to I_{W/X} \to 0$$ we have,
since $d>>0$, $\HH^1(T_Y\tensor
K_Y\tensor\mathcal{I}_{W/Y}(2d))\isom\HH^1(T_Y\tensor K_Y\tensor
I_{W/X}(2d)).$ We shall prove that the latter vanishes. For that we
use the exact sequence
$$0\to \Oh_{X}(-2)\to \Oh_X(-1)^{\oplus 2} \to I_{D/X} \to 0.$$
Tensoring this with $I_{C/X}(2d)$ yields
\comment{$T_Y\tensor K_Y\tensor I_{C/X}(2d)$ yields 
$$0\to T_Y\tensor K_Y\tensor I_{C/X}(2d-2)\to T_Y\tensor K_Y\tensor
I_{C/X}(2d-1)^{\oplus 2} \to T_Y\tensor K_Y\tensor I_{W/X}(2d) \to 0.$$}
$$0\to I_{C/X}(2d-2)\to I_{C/X}(2d-1)^{\oplus 2} \to I_{W/X}(2d) \to
0.$$ Since $\HH^j(T_Y\tensor K_Y\tensor{I}_{C/X}(2d-j))=0$ for
$j=1,~2$ from Lemma \ref{vanishings} above, we are done.
\end{proof}

\subsection{An auxiliary vector space} 
In this subsection, we shall construct an auxiliary vector space which surjects 
onto the domain of the map $\chi$. This construction is a crucial refinement of condition (1) in \cite{Wu1}
which was first proved in \cite{R}.
\begin{lemma}\label{utov}
Let $U:=\ker[\HH^0( T_{Y}\tensor K_{Y}(2d)) \to
\HH^0(K_{Y}(3d)_{|W})]$ and $V:=\ker[\HH^1(\Omega^1_{X|W}) \to
\HH^1(\Omega^2_{Y}(d)_{|W})]$. Then the natural map $U \to V$ is
a surjection.
\end{lemma}

\begin{proof}
Tensoring the short exact sequence
$$ 0 \to  T_X \to  T_{{Y|X}} \to \Oh_X(d) \to 0$$ with
$K_{Y}(2d)_{|W}$ and taking cohomology, we get
$$0 \to \HH^0( T_{X}\tensor K_{Y}(2d)_{|W}) \to
\HH^0( T_{Y}\tensor K_{Y}(2d)_{|W}) \to
\HH^0(K_{Y}(3d)_{|W}) $$

Since $ T_X\tensor K_{Y}(d) \isom \Omega^2_X$, we have
the following commutative diagram:
\begin{equation}\label{unexplained}
\begin{array}{ccccccc}
0 & \to & \UU  & \to & \HH^0( T_{Y}\tensor
K_{Y}(2d))
& \to & \HH^0(K_{Y}(3d)_{|W}) \\
& & \downarrow & & \downarrow & & || \\
0 & \to & \HH^0(\Omega^2_{X}(d)_{|W}) & \to &
 \HH^0( T_{Y}\tensor K_{Y}(2d)_{|W}) & \to &
\HH^0(K_{Y}(3d)_{|W}).\\
\end{array}
\end{equation}

The middle vertical arrow can be seen to be a surjection by using the
fact that the cokernel of this map injects into $\HH^1( T_{Y}\tensor
K_Y\tensor\mathcal I_{W/Y}(2d))$ which vanishes by Lemma
\ref{tgtvanishing}. By the snake lemma, the first map is also a
surjection. Since
$$\im[\HH^0(\Omega^2_{X}(d)_{|W}) \to \HH^1(\Omega^1_{X|W})]
= \ker[\HH^1(\Omega^1_{X|W}) \to
\HH^1(\Omega^2_{Y}(d)_{|W})]=\VV,$$ we have a surjection
$ \UU \onto \HH^0(\Omega^2_{X}(d)_{|W}) \onto \VV.$
\end{proof}

\subsection{The final lifting}
All that remains to be done now is to lift the elements from the auxiliary 
vector space $U$ constructed above to $\ker\gamma$ for which we need
the following
\begin{lemma}\label{multiplicationmap}
With notation as above, the multiplication map
$$\HH^0(\mathcal{I}_{W/Y}\tensor K_Y(2d))\tensor\HH^0(\Oh_{Y}(d)) \to
\HH^0(\mathcal{I}_{W/Y}\tensor K_Y(3d)) $$ is surjective.
\end{lemma}

\begin{proof}
Tensoring the exact sequence 
\begin{equation}\label{resforD}
0 \to \Oh_X(-2) \to \Oh_X(-1)^{\oplus 2} \to I_{D/X} \to 0 \, ,
\end{equation}
by $E$, we have
\begin{equation}\label{aa}
0 \to E(-2) \to E(-1)^{\oplus 2} \to I_{D/X}E \to 0 \, .
\end{equation}

Let $T_m:=\HH^0(\Oh_X(m))$. The exact sequence above gives rise to a diagram
with exact rows where the vertical arrows are all multiplication maps:
\[
\begin{array}{ccccccccc}
0 & \to &\HH^0(E(k-2))\tensor T_m & \to & 
\HH^0(E(k-1))^{\oplus{2}}\tensor T_m & \to  &
\HH^0(I_{D/X}E(k))\tensor T_m & \to & 0 \\
& &   \downarrow  &  & \downarrow &  & \downarrow  & & \\
0 & \to & \HH^0(E(m+k-2)) & \to & 
\HH^0(E(m+k-1))^{\oplus{2}} & \to &
\HH^0(I_{D/X}E(m+k)) & \to & 0. \\
\end{array}
\]
Since $E$ is $(\sum_{i=1}^{n-3}d_i-\alpha -n+3)$-regular by Lemma
\ref{c1inequality}, the middle vertical arrow is a surjection for
$k\geq (\sum_{i=1}^{n-3}d_i-\alpha-n+4)$ and $m\geq 0$. It follows
that the multiplication map
$$\HH^0(I_{D/X}E(k))\tensor\HH^0(\Oh_X(m))\to \HH^0(I_{D/X}E(m+k))$$
is surjective for $k\geq (\sum_{i=1}^{n-3}d_i-\alpha-n+4)$ and $m\geq
0$. Next consider the exact sequence $0 \to I_{D/X} \to I_{D/X}E \to
I_{W/X}(\alpha) \to 0$ obtained by tensoring sequence (\ref{serre}) by
$I_{D/X}(\alpha)$. Repeating the previous argument, it is easy to
check that the multiplication map
$$\HH^0(I_{W/X}(k))\tensor\HH^0(\Oh_X(m))\to \HH^0(I_{W/X}(m+k))$$ is
surjective for $k\geq (\sum_{i=1}^{n-3}d_i-n+4)$ and $m\geq 0$. Now
using the exact sequence
$$ 0 \to \Oh_{Y}(-d) \to \mathcal{I}_{W/Y} \to I_{W/X} \to 0,$$ and
the fact that the regularity of $\Oh_Y$ is $\sum_{i=1}^{n-4}d_i-n-2$,
we can conclude by repeating the argument above, that the
multiplication map
$$\HH^0(\mathcal{I}_{W/Y}\tensor K_Y(2d))\tensor\HH^0(\Oh_{Y}(d)) \to
\HH^0(\mathcal{I}_{W/Y}\tensor K_Y(3d))$$ is surjective.
\end{proof}

\section{Proofs of the main results}

Assume that a general, smooth complete intersection threefold $X\subset \bbP^{n}$ 
of sufficiently high multi-degree $(d_{1},\cdots, d_{n-3})$ 
supports an indecomposable ACM rank two vector
bundle $E$ with first Chern class $\alpha< d_{n-3}-1$.  This implies
that there exists a rank two bundle $\mathcal{E}$ on the universal
hypersurface $\mathcal{X} \subset Y \times S'$ where $S'$ is a Zariski
open subset of $S$, the moduli space of smooth, degree $d$
hypersurfaces of $Y$, such that for a general point $s\in S'$,
$\mathcal{E}_{|X_s}$ is normalised, indecomposable, ACM and its first
Chern class $\alpha_s=\alpha$ satisfies the inequality
above. Furthermore, from the construction of this family, one sees
that there exists a family of curves $\mathcal{C} \to S'$ such that
$\mathcal{C}_s$ is the zero locus of a section of
$\mathcal{E}_{|X_s}$. Let $\mathcal{Z}$ be a family of $1$-cycles with
fibre $\mathcal{Z}_s:=d\mathcal{C}_s-lD_s$ where, as before, $D_s$ is
plane section of $X_s$ and $l=l(s)$ is the degree of
$\mathcal{C}_s$.

\begin{prop}\label{finale}
In the situation above, $\delta\nu_{\mathcal{Z}} \not\equiv 0$.
\end{prop}

\begin{proof}
We shall show that $\delta\nu_{\mathcal{Z}}(s)\neq 0$ at any point
$s\in S$ parametrising a smooth hypersurface $X\subset Y$. To do
this, we shall lift elements of $U$ to $\ker\gamma$ in diagram
\ref{formula}. Since we have surjections $U \onto V \onto \bbC$, we
will be done. 

Let $\partial_f:\Omega^3_{Y}(2d)\to K_{Y}(3d)$ be the
derivative map where $f$ is the degree $d$ polynomial defining
$X$. Composing with the quotient $ K_{Y}(3d)/ K_{Y}(2d)$, we
get a map $\bar\partial_f:\Omega^3_{Y}(2d)\to
K_{Y}(3d)/K_{Y}(2d)$. Using the identification
$\Omega^3_{Y}\isom T_{Y}\tensor K_{Y}$, and taking
cohomology, we get
$$\HH^0( T_{Y}\tensor K_{Y}(2d))\by{\partial_f}
\HH^0(K_{Y}(3d)) \to
\frac{\HH^0(K_{Y}(3d))}{\HH^0(K_{Y}(2d))}.$$ The cokernel of
the composite map above can be identified with $\HH^2(\Omega^1_X)$
(see \cite{CGGH}, Page 174 or \cite{JL}, Chapter 9 for details).

The key ingredient in this lifting is the following commutative
diagram:
\begin{equation}
\begin{array}{ccc}
\HH^0(\Oh_{Y}(d))\tensor\HH^0(\mathcal I_{W/Y}\tensor
K_{Y}(2d)) & \stackrel{\gamma'}\to &
\frac{\HH^0(K_{Y}(3d))}{\partial_f\HH^0( T_{Y}\tensor
K_{Y}(2d))}\\ \downarrow & & \downarrow \\
\HH^1( T_{X})\tensor
\HH^1(\mathcal{I}_{W/Y}\tensor\Omega^2_{X}) & \by{\gamma} &
\HH^2(\Omega^1_{X}).\\
\end{array}
\end{equation}
Here the right vertical map is the one explained above. The horizontal
maps $\gamma$ and $\gamma'$ are (essentially) cup product maps. The
vertical map on the left is a tensor product of two maps: The first
factor is the composite $\HH^0(\Oh_{Y}(d))\to
\HH^0(\Oh_X(d))\to\HH^1( T_X)$. The normal bundle of $X\subset
Y$ is $\Oh_X(d)$ and $\HH^0(\Oh_X(d))\to \HH^1( T_X)$ is the
natural coboundary map in the cohomology sequence of the tangent
bundle sequence for this inclusion. The second factor $\HH^0(\mathcal
I_{W/Y}\tensor K_{Y}(2d))\to
\HH^1(\mathcal{I}_{W/Y}\tensor\Omega^2_{X})$ is also obtained
as above by observing that $T_X\tensor K_{Y}(d)\isom \Omega^2_X$.

This diagram yields a map $\ker\gamma' \to \ker\gamma$. To complete
the lifting, recall that by Lemma \ref{multiplicationmap}, the map
$\HH^0(\Oh_{Y}(d))\tensor\HH^0(\mathcal I_{W/Y}\tensor K_{Y}(2d))\onto
\HH^0(\mathcal{I}_{W/Y}\tensor{K_{Y}}(3d))$ is a surjection.
Restricting this map to $\ker{\gamma'}$, we get a surjection
$$\ker{\gamma'}\onto \bar{U}:=\partial_f\HH^0( T_{Y}\tensor
K_{Y}(2d))\cap \HH^0(\mathcal{I}_{W/Y}\tensor
K_{Y}(3d)).$$ 

Let $\widetilde{U}$ be the kernel of the map
$\HH^0( T_{Y}\tensor K_{Y}(2d)) \to
\HH^0(K_{Y}(3d)_{|X})$. Looking at the diagram analogous to
(\ref{unexplained}) obtained by replacing $W$ by $X$, we see that
there is a map $\widetilde{U}\to \HH^0(\Omega^2_X(d))$. The boundary
map $\HH^0(\Omega^2_X(d))\to \HH^1(\Omega^1_X)$ in the cohomology
sequence associated to diagram (\ref{2ndext}) is the zero map (this is
because the composite map $\HH^1(\Omega^1_X) \to
\HH^1(\Omega^2_{Y}(d)_{|X})\isom \HH^2(\Omega^2_{Y})$ is the
Gysin inclusion). This implies that the surjection $U \onto V$ of
Lemma \ref{utov} factors as $U \onto \bar{U}\onto U/\widetilde{U}\onto V$
and thus we have surjections $\ker\gamma'\onto \bar{U} \onto V
\stackrel{\chi}\onto \bbC$. By the compatibility of these maps with
the map $\ker\gamma'\to\ker\gamma$ and those in diagram
(\ref{formula}), we conclude (using Griffiths' formula) that
$\delta\nu_{\mathcal{Z}}(s)\neq 0$.
\end{proof}

\begin{proof}[Proof of Theorem B]
Assume that a general complete intersection threefold $X$ supports an
indecomposable normalised ACM bundle $E$, with $\alpha < d-1$. Let
$\mathcal{Z}$ be the family of degree zero $1$-cycles defined
earlier. By the refined Wu's criterion
$\delta\nu_{\mathcal{Z}}\not\equiv 0$: this contradicts the theorem of
Green and M\"uller-Stach. Thus we are done.
\end{proof}

\begin{proof}[Proof of Theorem A]
Let $X$ be a complete intersection subvariety of dimension four. Let
$E$ be an ACM bundle of rank two on $X$. As mentioned above, we may
assume $E$ to be normalised with first Chern class $\alpha$. Now
choose a general hypersurface $T\subset X$ of degree $d>>0$,
satisfying $d > \alpha +1$. Since $E\tensor\Oh_T$ is ACM, and $\alpha
< d-1$, it follows from Theorem B that $E\tensor\Oh_T$
splits. This implies by a standard argument that $E$ itself
splits. The case for dimension greater than four now follows in a 
similar way.
\end{proof}

\begin{proof}[Proof of Corollary C]
The proof of the first part is trivial. The proof of the second part
is as follows: Suppose $C=X\cap{\tilde{S}}$ where
$\tilde{S}\subset\bbP^n$ is a codimension two subscheme. Then
$C=X\cap{S}$ where $S:=\tilde{S}\cap{Y}$ is a surface in $Y$. We have
a commutative diagram

\[
\begin{array}{ccccccccc}
0&\to&\Omega^2_Y & \to&\Omega^2_Y(d)&\to&\Omega^2_Y(d)_{|X}&\to & 0 \\
 &   & \downarrow & & \downarrow & & \downarrow & & \\
0&\to&\Omega^2_{Y}\, _{|S} & \to&\Omega^2_Y(d)_{|S} 
&\to&\Omega^2_Y(d)_{{|C}}&\to & 0 \\
\end{array}
\]
Taking cohomology, we get a commutative diagram
\[
\begin{array}{ccc}
\HH^1(\Omega^2_Y(d)_{|X})& \isom  & \HH^2(\Omega^2_Y)\\
\downarrow &  & \downarrow \\ 
\HH^1(\Omega^2_Y(d)_{{|C}}) & \to & \HH^2(\Omega^2_{Y}\, _{|S}) \\
\end{array}
\]
The map $\HH^2(\Omega^2_Y) \to \HH^2(\Omega^2_Y\,_{|S})$ is non-zero
since the composite
$$\HH^2(\Omega^2_Y) \to \HH^2(\Omega^2_Y\,_{|S})\to
\HH^2(\Omega^2_S)$$ is a surjection which sends $h_Y^2 \mapsto h_S^2$
where $h_Y$ and $h_S$ are the classes of hyperplane sections in $Y$
and $S$ respectively. Thus the image of $h_X$ under the composite map
$$\HH^1(\Omega^1_X) \to \HH^1(\Omega^2_Y(d)_{|X}) \to
\HH^2(\Omega^2_{Y}\, _{|S})$$ is non-zero and hence its image under the
map 
$$\HH^1(\Omega^1_X) \to \HH^1(\Omega^2_Y(d)_{|X}) \to
\HH^1(\Omega^2_Y(d)_{|C})$$ is also non-zero. By the proof of Lemma
\ref{chiisonto}, if $E$ were indecomposable, then the above map is
zero. This implies when $\alpha<d-1$, that the associated rank two bundle splits, hence
$C$ is a complete intersection.
\end{proof}

\end{document}